\documentclass[a4paper,10pt]{article}
\usepackage{stmaryrd}
\usepackage{amsfonts}
\usepackage{bbm}
\usepackage{amscd}
\usepackage{mathrsfs}
\usepackage{latexsym,amssymb,amsmath,amscd,amscd,amsthm,amsxtra,xypic}
\usepackage[dvips]{graphicx}
\usepackage[utf8]{inputenc}
\usepackage[T1]{fontenc}
\usepackage{lmodern}
\usepackage{amssymb}
\usepackage[all]{xy}
\usepackage{nicefrac,mathtools,enumitem}
\usepackage{microtype}
\usepackage{lmodern}

\newcommand {\emptycomment}[1]{}
\newtheorem{thm}{Theorem}[section]

\newtheorem{pro}[thm]{Proposition}

\newtheorem{rmk}[thm]{Remark}
\newtheorem{defi}[thm]{Definition}

\newcommand{\lon }{\,\rightarrow\,}
\newcommand{\be }{\begin{equation}}
\newcommand{\ee }{\end{equation}}

\newcommand{\pf}{\noindent{\bf Proof.}\ }

\newcommand{\g}{\frkg}

\newcommand{\frkg}{\mathfrak g}
\newcommand{\frkh}{\mathfrak h}

\newcommand{\frkk}{\mathfrak k}
\newcommand{\frkl}{\mathfrak l}

\def\qed{\hfill ~\vrule height6pt width6pt depth0pt}

\newcommand{\br}[1]{   [ \cdot,    \cdot  ]   }

\newcommand{\gl}{\mathfrak {gl}}

\newcommand{\sgn}{\mathrm{sgn}}

\begin{document}
\title{
{ A complement on representations of Hom-Lie algebras } }
\author{  Xiong Zhen\\
Department of Mathematics and Computer, Yichun University,
 Jiangxi 336000,  China
}
\date{email:205137@jxycu.edu.cn}
\footnotetext{{\it{Keyword}: Hom-Lie algebras; representations; coboundary operators
 }}
 \footnotetext{{\it{MSC}}: 17B99, 55U15}
 \footnotetext{Supported by the National Natural Science Foundation of China (No.11771382) and  the Science and Technology Project(GJJ161029)of Department of Education, Jiangxi Province.}

\maketitle

\begin{abstract}
 In this paper, we give a new series of coboundary operators of Hom-Lie algebras. And prove that cohomology groups with respect to coboundary operators are isomorphic. Then, we revisit representations of Hom-Lie algebras, and prove that there is an one-to-one correspondence between Hom-Lie algebraic structure on vector space $\frkg$ and these coboundary operators on $\Lambda\frkg^*\otimes V$.
\end{abstract}

%\tableofcontents

\section{Introduction}

The notion of Hom-Lie algebras was introduced by Hartwig, Larsson,
and Silvestrov in \cite{HLS} as part of a study of deformations of
the Witt and the Virasoro algebras. In a Hom-Lie algebra, the Jacobi
identity is twisted by a linear map, called the Hom-Jacobi identity.
Some $q$-deformations of the Witt and the Virasoro algebras have the
structure of a Hom-Lie algebra \cite{HLS,hu}. Because of close relation
to discrete and deformed vector fields and differential calculus
\cite{HLS,LD1,LD2},   more people pay special attention to this algebraic structure.
In particular, Hom-Lie algebras on semisimple Lie algebras are studied in \cite{jin1};
Its geometric generalization is given in \cite{cl,hom-Lie algebroids}; Quadratic Hom-Lie algebras are studied in \cite{BM};
Representation theory, cohomology and homology theory are systematically studied in \cite{D,homlie1,chengsu,Yao2,AEM}; Bialgebra theory
and Hom-(Classical) Yang-Baxter Equation are studied in \cite{MS3,shengbai,Yao1,YaoCYBE}; The notion of a Hom-Lie 2-algebra,
which is a categorification of a Hom-Lie algebra, is introduced in \cite{homlie2}, in which the relation with Hom-left-symmetric algebras  \cite{MS2}
and  symplectic Hom-Lie algebras are studied.

Let $(\frkg,[\cdot,\cdot],\alpha)$ be a Hom-Lie algebra, $V$ be a
vector space, $\rho:\frkg\longrightarrow\frkg\frkl(V)$ be a
representation of $(\frkg,[\cdot,\cdot],\alpha)$ on the vector space
$V$ with respect to $\beta\in\frkg\frkl(V)$. The set of {\bf
$k$-cochains} on $\frkg$ with values in $V$, which we denote by
$C^k(\frkg;V)$, is the set of skewsymmetric $k$-linear maps from
$\frkg\times\cdots\times\frkg$($k$-times) to $V$:
$$C^k(\frkg;V):=\{\eta:\wedge^k\frkg\longrightarrow V ~\mbox{is a
linear map}\}.$$
In \cite{yx}, the authors define Hom-$k$-cochains: $C_{\alpha,\beta}^k(\frkg;V)=\{\eta\in C^k(\frkg;V)|\eta\circ\alpha=\beta\circ\eta\}$, which is a subset of $C^k(\frkg;V)$.  There are a series of coboundary operators $\hat{d}^s$ define on $C_{\alpha,\beta}^k(\frkg;V)$; in \cite{cy}, the authors give a special coboundary operator of regular Hom-Lie algebras. For regular Hom-Lie algebras, there are many works are done by the special coboundary operator \cite{cl,cy} .
In this article, we give a new series coboundary operators on $k$-cochains $C^k(\frkg;V)$, prove that cohomology groups with respect to these coboundary operators, are isomorphic. Then,
we revisit representations of Hom-Lie algebras, and
 generalize  the
result
"If $\frkk$ is a Lie algebra, $\rho:\frkk\longrightarrow\gl(V)$ is
a representation if and only if there is a degree-$1$ operator $D$
on $\Lambda\frkk^*\otimes V$ satisfying $D^2=0$, and
  $$
  D(\xi\wedge \eta\otimes u)=d_\frkk\xi\wedge\eta\otimes u+(-1)^k\xi \wedge D(\eta\otimes u),\quad \forall~\xi\in \wedge^k\frkk^*,~ \eta\in\wedge^l\frkk^*,~u\in V,
  $$
  where $d_\frkk:\wedge^k\g^*\longrightarrow \wedge^{k+1}\g^*$ is the coboundary operator associated to the trivial
  representation."

The paper is organized as follows. In Section 2, we first recall some
necessary background knowledge: Hom-Lie algebras and their
representations. Then, we show that $d^s$ is coboundary operators of Hom-Lie algebras and prove that cohomology groups with respect to these coboundary operators, are isomorphic(Theorem \ref{thm2}). In Section 3, we give some properties of $d^s$ and revisit representations of Hom-Lie algebras, and have Theorem \ref{thm3}.

\section{Cohomology operators of Hom-Lie algebras }

The notion of a Hom-Lie algebra was introduced in \cite{HLS}, see also \cite{BM,MS2} for more information.
\begin{defi}
\begin{itemize}
\item[\rm(1)]
  A Hom-Lie algebra is a triple $(\frkg,\br ,,\alpha)$ consisting of a
  vector space $\frkg$, a skewsymmetric bilinear map (bracket) $\br,:\wedge^2\frkg\longrightarrow
  \frkg$ and a linear transformation $\alpha:\frkg\lon\frkg$ satisfying $\alpha[x,y]=[\alpha(x),\alpha(y)]$, and the following Hom-Jacobi
  identity:
  \begin{equation}
   [\alpha(x),[y,z]]+[\alpha(y),[z,x]]+[\alpha(z),[x,y]]=0,\quad\forall
x,y,z\in\frkg.
  \end{equation}

 A Hom-Lie algebra is called a regular Hom-Lie algebra if $\alpha$ is
a linear automorphism.

 \item[\rm(2)] A subspace $\frkh\subset\frkg$ is a Hom-Lie sub-algebra of $(\frkg,\br ,,\alpha)$ if
 $\alpha(\frkh)\subset\frkh$ and
  $\frkh$ is closed under the bracket operation $\br,$, i.e. for all $ x,y\in\frkh$,
  $[x,y] \in\frkh.  $
  \item[\rm(3)] A morphism from the  Hom-Lie algebra
$(\frkg,[\cdot,\cdot]_{\frkg},\alpha)$ to the hom-Lie algebra
$(\frkh,[\cdot,\cdot]_{\frkh},\gamma)$ is a linear map
$\psi:\frkg\longrightarrow\frkh$ such that
$\psi([x,y]_{\frkg})=[\psi(x),\psi(y)]_{\frkh}$ and
$\psi\circ \alpha =\gamma\circ \psi$.
  \end{itemize}
\end{defi}

 Representation and cohomology theories of Hom-Lie algebra are
systematically introduced in \cite{AEM,homlie1}. See \cite{Yao2} for
homology theories of Hom-Lie algebras.

\begin{defi}
  A representation of the Hom-Lie algebra $(\frkg,\br,,\alpha)$ on
  a vector space $V$ with respect to $\beta\in\gl(V)$ is a linear map
  $\rho:\frkg\longrightarrow \gl(V)$, such that for all
  $x,y\in\frkg$, the following equalities are satisfied:
  \begin{eqnarray}
 \label{representation1} \rho (\alpha(x))\circ \beta&=&\beta\circ \rho (x);\\\label{representation2}
    \rho([x,y] )\circ
    \beta&=&\rho (\alpha(x))\circ\rho (y)-\rho (\alpha(y))\circ\rho (x).
  \end{eqnarray}
\end{defi}

Let $(\frkg,[\cdot,\cdot],\alpha)$ be a Hom-Lie algebra, $V$ be a
vector space, $\rho:\frkg\longrightarrow\frkg\frkl(V)$ be a
representation of $(\frkg,[\cdot,\cdot],\alpha)$ on the vector space
$V$ with respect to $\beta\in GL(V)$, where $\beta$ is invertible. In this paper, we just consider $\beta\in GL(V)$.

The set of {\bf
$k$-cochains} on $\frkg$ with values in $V$, which we denote by
$C^k(\frkg;V)$, is the set of skewsymmetric $k$-linear maps from
$\frkg\times\cdots\times\frkg$($k$-times) to $V$:
$$C^k(\frkg;V):=\{\eta:\wedge^k\frkg\longrightarrow V ~\mbox{is a
linear map}\}.$$
For $s=0,1,2,\ldots,$ define
$d^s:C^k(\frkg;V)\longrightarrow
C^{k+1}(\frkg;V)$ by
\begin{eqnarray*}
d^s\eta(x_1,\cdots,x_{k+1})&=&\sum_{i=1}^{k+1}(-1)^{i+1}\beta^{k+1+s}\rho(x_i)\beta^{-k-2-s}\eta(\alpha(x_1),\cdots,\hat{x_i},\cdots,\alpha(x_{k+1}))\\
&&+\sum_{i<j}(-1)^{i+j}\eta([x_i,x_j],\alpha(x_1),\cdots,\widehat{x_{i,j}},\cdots,\alpha(x_{k+1})),
\end{eqnarray*}
where $\beta^{-1}$ is the inverse of $\beta$, $\eta\in C^k(\frkg;V)$.
\begin{pro}
With the above notations, the map $d^s$ is a coboundary operator, i.e. $d^s\circ d^s=0$.
\end{pro}
\pf For any $\eta\in C^k(\frkg;V)$, by straightforward computations, we have
\begin{eqnarray*}
d^s\circ d^s\eta(x_1,\cdots,x_{k+2})&=&\sum_{i=1}^{k+2}(-1)^{i+1}\beta^{k+2+s}\rho(x_i)\beta^{-k-3-s}d^s\eta(\alpha(x_1),\cdots,\hat{x_i},\cdots,\alpha(x_{k+2}))\\
&&+\sum_{i<j}(-1)^{i+j}d^s\eta([x_i,x_j],\alpha(x_1),\cdots,\widehat{x_{i,j}},\cdots,\alpha(x_{k+2})).
\end{eqnarray*}
And\\
$d^s\eta(\alpha(x_1),\cdots,\hat{x_i},\cdots,\alpha(x_{k+2}))$
\begin{eqnarray*}
&=&\sum_{l<i}(-1)^{l+1}\beta^{k+1+s}\rho(\alpha(x_l))\beta^{-k-2-s}\eta(\alpha^2(x_1),\cdots,\widehat{x_{l,i}},\cdots,\alpha^2(x_{k+2}))\\
&&+\sum_{l>i}(-1)^{l}\beta^{k+1+s}\rho(\alpha(x_l))\beta^{-k-2-s}\eta(\alpha^2(x_1),\cdots,\widehat{x_{i,l}},\cdots,\alpha^2(x_{k+2}))\\
&&+\sum_{m<n<i}(-1)^{m+n}\eta(\alpha([x_i,x_j]),\alpha^2(x_1),\cdots,\widehat{x_{m,n,i}},\cdots,\alpha^2(x_{k+2}))\\
&&+\sum_{m<i<n}(-1)^{m+n-1}\eta(\alpha([x_i,x_j]),\alpha^2(x_1),\cdots,\widehat{x_{m,i,n}},\cdots,\alpha^2(x_{k+2}))\\
&&+\sum_{i<m<n}(-1)^{m+n}\eta(\alpha([x_i,x_j]),\alpha^2(x_1),\cdots,\widehat{x_{i,m,n}},\cdots,\alpha^2(x_{k+2})).
\end{eqnarray*}
At the same time, we have\\
$d^s\eta([x_i,x_j],\alpha(x_1),\cdots,\widehat{x_{i,j}},\cdots,\alpha(x_{k+2})$
\begin{eqnarray}
&=&\beta^{k+1+s}\rho([x_i,x_j])\beta^{-k-2-s}\eta(\alpha^2(x_1),\cdots,\widehat{x_{i,j}},\cdots,\alpha^2(x_{k+2}))\nonumber\\
&&+\sum_{p<i<j}(-1)^{p}\beta^{k+1+s}\rho(\alpha(x_p))\beta^{-k-2-s}\eta(\alpha([x_i,x_j]),\alpha^2(x_1),\cdots,\widehat{x_{p,i,j}},\cdots,\alpha^2(x_{k+2}))\nonumber\\
&&+\sum_{i<p<j}(-1)^{p+1}\beta^{k+1+s}\rho(\alpha(x_p))\beta^{-k-2-s}\eta(\alpha([x_i,x_j]),\alpha^2(x_1),\cdots,\widehat{x_{i,p,j}},\cdots,\alpha^2(x_{k+2}))\nonumber\\
&&+\sum_{i<j<p}(-1)^{p}\beta^{k+1+s}\rho(\alpha(x_p))\beta^{-k-2-s}\eta(\alpha([x_i,x_j]),\alpha^2(x_1),\cdots,\widehat{x_{i,j,p}},\cdots,\alpha^2(x_{k+2}))\nonumber\\
&&+\sum_{q<i<j}(-1)^{1+q}\eta([[x_i,x_j],\alpha(x_q)],\alpha^2(x_1),\cdots,\widehat{x_{q,i,j}},\cdots,\alpha^2(x_{k+2}))\label{eq1}\\
&&+\sum_{i<q<j}(-1)^{q}\eta([[x_i,x_j],\alpha(x_q)],\alpha^2(x_1),\cdots,\widehat{x_{i,q,j}},\cdots,\alpha^2(x_{k+2}))\label{eq2}\\
&&+\sum_{i<j<q}(-1)^{1+q}\eta([[x_i,x_j],\alpha(x_q)],\alpha^2(x_1),\cdots,\widehat{x_{i,j,q}},\cdots,\alpha^2(x_{k+2}))\label{eq3}\\
&&+\sum_{m<n<i<j}(-1)^{m+n}\eta([\alpha(x_m),\alpha(x_n)],\alpha([x_i,x_j]),\alpha^2(x_1),\cdots,\widehat{x_{m,n,i,j}},\cdots,\alpha^2(x_{k+2}))\label{eq4}\\
&&+\cdots\nonumber
\end{eqnarray}
By Hom-Jacobi  identity:
$$(\ref{eq1})+(\ref{eq2})+(\ref{eq3})=0,$$
and we have: $(\ref{eq4})+\cdots=0.$

By $\rho(\alpha(x))\beta=\beta\rho(x)$ and $\rho(\alpha(x))=\beta\rho(x)\beta^{-1}$, we have\\
$d^s\circ d^s\eta(x_1,\cdots,x_{k+2})$
\begin{eqnarray}
&=&\sum_{l<i}(-1)^{l+i}\beta^{k-1+s}\rho(\alpha^3(x_i))\rho(\alpha^2(x_l))\beta^{-k-1-s}\eta(\alpha^2(x_1),\cdots,\widehat{x_{l,i}},\cdots,\alpha^2(x_{k+2}))\label{eq5}\\
&&+\sum_{l>i}(-1)^{l+i+1}\beta^{k-1+s}\rho(\alpha^3(x_i))\rho(\alpha^2(x_l))\beta^{-k-1-s}\eta(\alpha^2(x_1),\cdots,\widehat{x_{i,l}},\cdots,\alpha^2(x_{k+2}))\label{eq6}\\
&&+\sum_{m<n<i}(-1)^{m+n+i+1}\beta^{k+2+s}\rho(x_i)\beta^{-k-3-s}\eta([\alpha(x_i),\alpha(x_j)],\alpha^2(x_1),\cdots,\widehat{x_{m,n,i}},\cdots,\alpha^2(x_{k+2}))\nonumber\\
&&+\sum_{m<i<n}(-1)^{m+n+i}\beta^{k+2+s}\rho(x_i)\beta^{-k-3-s}\eta([\alpha(x_i),\alpha(x_j)],\alpha^2(x_1),\cdots,\widehat{x_{m,i,n}},\cdots,\alpha^2(x_{k+2}))\nonumber\\
&&+\sum_{i<m<n}(-1)^{m+n+i+1}\beta^{k+2+s}\rho(x_i)\beta^{-k-3-s}\eta([\alpha(x_i),\alpha(x_j)],\alpha^2(x_1),\cdots,\widehat{x_{i,m,n}},\cdots,\alpha^2(x_{k+2}))\nonumber\\
&&+\sum_{i<j}(-1)^{i+j}\beta^{k-1+s}\rho([\alpha^2(x_i),\alpha^2(x_j)])\beta\beta^{-k-1-s}\eta(\alpha^2(x_1),\cdots,\widehat{x_{i,j}},\cdots,\alpha^2(x_{k+2}))\label{eq7}\\
&&+\sum_{p<i<j}(-1)^{p+i+j}\beta^{k+2+s}\rho(x_p)\beta^{-k-3-s}\eta([\alpha(x_i),\alpha(x_j)],\alpha^2(x_1),\cdots,\widehat{x_{p,i,j}},\cdots,\alpha^2(x_{k+2}))\nonumber\\
&&+\sum_{i<p<j}(-1)^{p+i+j+1}\beta^{k+2+s}\rho(x_p)\beta^{-k-3-s}\eta([\alpha(x_i),\alpha(x_j)],\alpha^2(x_1),\cdots,\widehat{x_{i,p,j}},\cdots,\alpha^2(x_{k+2}))\nonumber\\
&&+\sum_{i<j<p}(-1)^{p+i+j}\beta^{k+2+s}\rho(x_p)\beta^{-k-3-s}\eta([\alpha(x_i),\alpha(x_j)],\alpha^2(x_1),\cdots,\widehat{x_{i,j,p}},\cdots,\alpha^2(x_{k+2}))\nonumber.
\end{eqnarray}
By $\rho([x,y])\beta=\rho(\alpha(x))\rho(y)-\rho(\alpha(y))\rho(x)$, we have
$$(\ref{eq5})+(\ref{eq6})+(\ref{eq7})=0.$$
About above equations, sum of the rest six equations is zero. So, we proof that $d^s\circ d^s=0$.\qed

\begin{rmk}
In \cite{yx}, the authors give a series coboundary operators $\hat{d}^s:C_{\alpha,\beta}^k(\frkg;V)\longrightarrow C_{\alpha,\beta}^{k+1}(\frkg;V)$,
on the set of Hom-$k$-cochains $C_{\alpha,\beta}^k(\frkg;V)=\{\eta\in C^k(\frkg;V)|\eta\circ\alpha=\beta\circ\eta\}$. The coboundary operators $d^s$ we define are not the same as those in \cite{cy}, \cite{cy} just give coboundary operator for regular Hom-Lie algebras.
\end{rmk}

From $d^s$, we know that the coboundary operator associated to the trivial representation is $d:\wedge^k\frkg^*\longrightarrow\wedge^{k+1}\frkg^*$,
$$d\xi(x_1,\cdots,x_{k+1})=\sum_{i<j}(-1)^{i+j}\xi([x_i,x_j],\alpha(x_1),\cdots,\widehat{x_{i,j}},\cdots,\alpha(x_{k+1})).$$
For Hom-Lie algebra $(\frkg,[\cdot,\cdot],\alpha)$ and representation $\rho$ with respect with $\beta$, $\alpha$ induces a map
$\bar{\alpha}:C^l(\frkg;V)\longrightarrow C^{l}(\frkg;V)$ via
$$\bar{\alpha}(\eta)(x_1,\cdots,x_l)=\eta(\alpha(x_1),\cdots,\alpha(x_l)).$$
And $\beta$ induces a map $\bar{\beta}:C^l(\frkg;V)\longrightarrow C^{l}(\frkg;V)$ via
$$\bar{\beta}(\eta)(x_1,\cdots,x_l)=\beta\circ\eta(x_1,\cdots,x_l).$$
$C^\bullet(\frkg;V)=\oplus_lC^l(\frkg;V)$ is a $\wedge^\bullet=\oplus_k\wedge^k\frkg^*$-module, where the action $\diamond:\wedge^k\frkg^*\times
C^l(\frkg;V)\longrightarrow
C^{k+l}(\frkg;V)$ is given by
 $$
 \xi\diamond\eta(x_1,\cdots,x_{k+l})=\sum_\kappa \sgn(\kappa)\eta(x_{\kappa(1)},\cdots,x_{\kappa(k)})\eta(x_{\kappa(k+1)},\cdots,x_{\kappa(k+l)}),
 $$
where $\xi\in \wedge^k\frkg^*, \eta\in
C^l(\frkg;V),$ and the summation is taken over  $(k,l)$-unshuffles.

Obviously, for $\xi,\xi_1,\xi_2\in\wedge^\bullet, \eta\in C^\bullet(\frkg;V)$, we have
\begin{eqnarray*}
\bar{\alpha}(\xi_1\wedge\xi_2)&=&\bar{\alpha}(\xi_1)\wedge\bar{\alpha}(\xi_2);\\
\bar{\alpha}(\xi\diamond\eta)&=&\bar{\alpha}(\xi)\diamond\bar{\alpha}(\eta);\\
\bar{\beta}(\xi\diamond\eta)&=&\xi\diamond\bar{\beta}(\eta).
\end{eqnarray*}

Associated to the representation $\rho$, we obtain the complex$(C^k(\frkg:V),d^s)$. Denote the set
of closed $k$-cochains by $Z^k(d^s)$ and the set of exact $k$-cochains by $B^k(d^s)$. Denote
the corresponding cohomology by
$$H^k(d^s)=Z^k(d^s)/B^k(d^s).$$
Now, we study the relation between $H^k(d^s)$ and $H^k(d^{s+1})$.
\begin{pro}
With the above notations, we have
$$\bar{\beta}\circ d^s=d^{s+1}\circ\bar{\beta}.$$
\end{pro}
\pf For $\eta\in C^l(\frkg;V)$, we have
\begin{eqnarray*}
\bar{\beta}\circ d^s\eta(x_1,\cdots,x_{l+1})&=&\sum_{i=1}^{l+1}(-1)^{i+1}\beta^{l+2+s}\rho(x_i)\beta^{-l-2-s}\eta(\alpha(x_1),\cdots,\hat{x_i},\cdots,\alpha(x_{l+1}))\\
&&+\sum_{i<j}(-1)^{i+j}\beta\circ\eta([x_i,x_j],\alpha(x_1),\cdots,\widehat{x_{i,j}},\cdots,\alpha(x_{l+1}))\\
&=&\sum_{i=1}^{l+1}(-1)^{i+1}\beta^{l+2+s}\rho(x_i)\beta^{-l-3-s}\bar{\beta}(\eta)(\alpha(x_1),\cdots,\hat{x_i},\cdots,\alpha(x_{l+1}))\\
&&+\sum_{i<j}(-1)^{i+j}\bar{\beta}(\eta)([x_i,x_j],\alpha(x_1),\cdots,\widehat{x_{i,j}},\cdots,\alpha(x_{l+1}))\\
&=&d^{s+1}(\bar{\beta}(\eta))(x_1,\cdots,x_{l+1}),
\end{eqnarray*}
which implies that $\bar{\beta}\circ d^s=d^{s+1}\circ\bar{\beta}.$\qed

\begin{thm}\label{thm2}
For $s=0,1,2\ldots$, we have: $H^k(d^s)\cong H^k(d^{s+1}).$
\end{thm}
\pf By $d^{s+1}\circ \bar{\beta}=\bar{\beta}\circ d^{s}$, for
$\eta\in Z^k(d^{s})$, we have $\bar{\beta}(\eta)\in Z^k(d^{s+1})$.
On the other hand, for $\eta_1\in B^k(d^{s})$, there is $\omega\in
C^{k-1}(\frkg;V)$, such that: $ \eta_1=d^{s}\omega$.
 so,
$$\bar{\beta}(\eta_1)=\bar{\beta}\circ d^{s}\omega=d^{s+1}\circ
\bar{\beta}(\omega).$$
Obviously, $\bar{\beta}(\omega)\in
C^{k-1}(\frkg;V)$, then, $\bar{\beta}(\eta_1)\in
B^k(d^{s+1})$. Actually, we have proof:
$$\bar{\beta}(Z^k(d^{s}))\subset Z^k(d^{s+1}),\quad
\bar{\beta}(B^k(d^{s}))\subset B^k(d^{s+1}).$$
Next, for $\beta^{-1}$, we define map $\overline{\beta^{-1}}:C^k(\frkg;V)\rightarrow
C^k(\frkg;V)$ by
$$\overline{\beta^{-1}}(\eta)(x_1,\cdots,x_k)=\beta^{-1}\circ
\eta(x_1,\ldots,x_k),\quad \forall \eta\in
C^k(\frkg;V).$$
For $\eta\in Z^k(d^{s+1})$, we have $d^{s+1}\eta=0$. By
$$\bar{\beta}\circ
d^{s}\circ
\overline{\beta^{-1}}(\eta)=d^{s+1}\circ\bar{\beta}\circ\overline{\beta^{-1}}(\eta)=d^{s+1}\eta=0.$$
We have:
 $$d^{s}\circ \overline{\beta^{-1}}(\eta)=0,$$
 then, we have:
$$\overline{\beta^{-1}}(\eta)\in Z^k(d^{s}).$$
 On the other hand, for $\eta_1\in
B^k(d^{s+1})$, there is $\omega\in
C^{k-1}(\frkg;V)$, such that $\eta_1=d^{s+1}\omega$.
Then, we have:
$$\bar{\beta}\circ d^{s}\circ \overline{\beta^{-1}}(\omega)=d^{s+1}\circ\bar{\beta}\circ\overline{\beta^{-1}}(\omega)=d^{s+1}\omega=\eta_1.$$
So, $$d^{s}\circ\overline{\beta^{-1}}(\omega)=\beta^{-1}\circ
\eta_1,$$ then $$\overline{\beta^{-1}}(\eta_1)\in B^k(d^{s}).$$
Actually, we have proof:
$$\overline{\beta^{-1}}(Z^k(d^{s+1}))\subset Z^k(d^{s});\quad \overline{\beta^{-1}}(B^k(d^{s+1}))\subset B^k(d^{s}).$$
Now, we complete the proof.\qed

\section{Representations of Hom-Lie algebras-revisited}

We first consider the coboundary operator $d$, which is associated to the trivial representation. The following is right.
\begin{pro}
For $\xi_1\in\wedge^k\frkg^*,\xi_2\in\wedge^l\frkg^*$, we have
$$d(\xi_1\wedge\xi_2)=d\xi_1\wedge\bar{\alpha}(\xi_2)+(-1)^k\bar{\alpha}(\xi_1)\wedge d\xi_2.$$
\end{pro}
\pf  This proof is similar to Proposition 3.2 in \cite{yx}.\qed

\begin{pro}
For $\xi\in\wedge^k\frkg^*, \eta\in C^l(\frkg;V)$, we have
$$d^s(\xi\diamond\eta)=d\xi\diamond\bar{\alpha}(\eta)+(-1)^k\bar{\alpha}(\eta)\diamond d^{s+k}\eta.$$
\end{pro}
\pf First let $k=1$, then $\xi\diamond\eta\in
C^{l+1}(\frkg;V)$. We have
\begin{eqnarray*}
&&d^s(\xi\diamond\eta)(x_1,\cdots,x_{l+2})\\
&=&\sum_{i=1}^{k+2}(-1)^{l+1}\beta^{l+2+s}\rho(x_i)\beta^{-l-3-s}\xi\diamond\eta(x_1,\cdots,\hat{x_i},\cdots,x_{l+2})\\
&&+\sum_{i<j}(-1)^{i+j}\xi\diamond\eta([x_i,x_j],\alpha(x_1),\cdots,\widehat{x_{i,j}},\cdots,\alpha(x_{l+2}))\\
&=&\sum_{i<j}(-1)^{i+j}\xi([x_i,x_j])\eta(\alpha(x_1),\cdots,\widehat{x_{i,j}},\cdots,\alpha(x_{l+2}))\\
&&+\sum_{q<i}(-1)^{i+q}\xi(\alpha(x_i))\beta^{l+2+s}\rho(x_q)\beta^{-l-3-s}\eta(\alpha(x_1),\cdots,\widehat{x_{q,i}},\cdots,\alpha(x_{l+2}))\\
&&+\sum_{i<q}(-1)^{i+q+1}\xi(\alpha(x_i))\beta^{l+2+s}\rho(x_q)\beta^{-l-3-s}\eta(\alpha(x_1),\cdots,\widehat{x_{i,q}},\cdots,\alpha(x_{l+2}))\\
&&+\sum_{q<i<j}(-1)^{q+i+j}\xi(\alpha(x_q))\eta([x_i,x_j],\alpha(x_1),\cdots,\widehat{x_{q,i,j}},\cdots,\alpha(x_{l+2}))\\
&&+\sum_{i<q<j}(-1)^{q+i+j+1}\xi(\alpha(x_q))\eta([x_i,x_j],\alpha(x_1),\cdots,\widehat{x_{i,q,j}},\cdots,\alpha(x_{l+2}))\\
&&+\sum_{i<j<q}(-1)^{q+i+j}\eta(\alpha(x_q))\eta([x_i,x_j],\alpha(x_1),\cdots,\widehat{x_{i,j,q}},\cdots,\alpha(x_{l+2}))\\
&=&d \xi\diamond
\bar{\alpha}(\eta)(x_1,\cdots,x_{l+2})+(-1)^1\bar{\alpha}(\xi)\diamond
d^{s+1}\varphi(x_1,\cdots,x_{l+2}).
\end{eqnarray*}

Thus, when $k=1$, we have
$$d^s(\xi\diamond\eta)=d \xi\diamond
\bar{\alpha}(\eta)+(-1)^1\bar{\alpha}(\xi)\diamond d^{s+1}\eta.$$
 By induction on $k$,  assume that when $k=n$, we have
$$d^s(\xi\diamond\eta)=d \xi\diamond\bar{\alpha}(\eta)+(-1)^n\bar{\alpha}(\xi)\diamond d^{s+n}\eta.$$
For $\omega\in \frkg^*$, $\xi\wedge\omega\in
 \wedge^{n+1}\frkg^*$,  we have
\begin{eqnarray*}
d^s((\xi\wedge\omega)\diamond\eta)&=&d^s(\xi\diamond(\omega\diamond\eta))\\
&=&d \xi\diamond\bar{\alpha}(\omega\diamond\eta)+(-1)^n\bar{\alpha}(\xi)\diamond d^{s+n}(\omega\diamond\eta)\\
&=&(d \xi\wedge\omega)\diamond\bar{\alpha}(\eta)+(-1)^n\bar{\alpha}(\xi)\diamond(d \omega\diamond\bar{\alpha}(\eta)+(-1)\bar{\alpha}(\omega)\diamond d^{s+n+1}\eta)\\
&=&(d \xi\wedge\bar{\alpha}(\omega)+(-1)^n\bar{\alpha}(\xi)\wedge
d \omega)\diamond\bar{\alpha}(\eta)+(-1)^{n+1}\bar{\alpha}(\xi\wedge\omega)\diamond
d^{s+n+1}\eta\\
&=&d (\eta\wedge\omega)\diamond\bar{\alpha}(\eta)+(-1)^{n+1}\bar{\alpha}(\eta\wedge\omega)\diamond
d^{s+n+1}\eta.
\end{eqnarray*}
The proof is completed.\qed

\begin{pro}
With the above notations, we have
$$\bar{\alpha}\circ d^s=d^{s+1}\circ\bar{\alpha}.$$
\end{pro}
\pf For any $\eta\in C^l(\frkg;V)$, by $\rho(\alpha(x_i))=\beta\circ\rho(x_i)\circ\beta$, we have
\begin{eqnarray*}
\bar{\alpha}\circ d^s\eta(x_1,\cdots,x_{l+1})&=&d^s\eta(\alpha(x_1),\cdots,\alpha(x_{l+1}))\\
&=&\sum_{i=1}^{l+1}(-1)^{i+1}\beta^{l+1+s}\rho(\alpha(x_i))\beta^{-l-2-s}\eta(\alpha^2(x_1),\cdots,\hat{x_i},\cdots,\alpha^2(x_{l+1}))\\
&&+\sum_{i<j}(-1)^{i+j}\eta([\alpha(x_i),\alpha(x_j)],\alpha^2(x_1),\cdots,\widehat{x_{i,j}},\cdots,\alpha^2(x_{l+1}))\\
&=&\sum_{i=1}^{l+1}(-1)^{i+1}\beta^{l+2+s}\rho(x_i)\beta^{-l-3-s}\bar{\alpha}(\eta)(\alpha(x_1),\cdots,\hat{x_i},\cdots,\alpha(x_{l+1}))\\
&&+\sum_{i<j}(-1)^{i+j}\bar{\alpha}(\eta)([x_i,x_j],\alpha(x_1),\cdots,\widehat{x_{i,j}},\cdots,\alpha(x_{l+1}))\\
&=&d^{s+1}(\bar{\alpha}(\eta))(x_1,\cdots,x_{l+1}),
\end{eqnarray*}
which implies that $\bar{\alpha}\circ d^s=d^{s+1}\circ\bar{\alpha}.$\qed

The converse of the above conclusions are also true. Thus, we have the following theorem, which generalize  the
result
" If $\frkk$ is a Lie algebra, $\rho:\frkk\longrightarrow\gl(V)$ is
a representation if and only if there is a degree-$1$ operator $D$
on $\Lambda\frkk^*\otimes V$ satisfying $D^2=0$, and
  $$
  D(\xi\wedge \eta\otimes u)=d_\frkk\xi\wedge\eta\otimes u+(-1)^k\xi \wedge D(\eta\otimes u),\quad \forall~\xi\in \wedge^k\frkk^*,~ \eta\in\wedge^l\frkk^*,~u\in V,
  $$
  where $d_\frkk:\wedge^k\g^*\longrightarrow \wedge^{k+1}\g^*$ is the coboundary operator associated to the trivial
  representation."
\begin{thm}\label{thm3}
Let $V$ be a vector space, $\beta\in GL(V)$. Then
$(\frkg,[\cdot,\cdot],\alpha)$ is a Hom-Lie algebra, and
$\rho:\frkg\longrightarrow\frkg\frkl(V)$ is a representation of
$(\frkg,[\cdot,\cdot],\alpha)$ on the vector space $V$ with respect
to $\beta$ if and only if there exists:
$d^s:C^l(\frkg;V)\longrightarrow C^{l+1}(\frkg;V)$, $s=0,1,2,\ldots$ and such
that:
\begin{itemize}
\item[\rm{i)}]$d^s\circ d^s=0$;
\item[\rm{ii)}]for any $\xi\in \wedge^k\frkg^*,\eta\in C^l(\frkg;V)$, we have
$$d^s(\xi\diamond\eta)=d\xi\diamond\bar{\alpha}(\eta)+(-1)^k\bar{\alpha}(\xi)\diamond d^{s+k}\eta;$$
where $d:\wedge^k\frkg^*\longrightarrow \wedge^{k+1}\frkg^*$ is the coboundary operator associated to the trivial
  representation.
\item[\rm{iii)}]$\bar{\alpha}\circ d^s=d^{s+1}\circ\bar{\alpha}$.
\end{itemize}
\end{thm}

\pf With Propositions we have proof above. We just need to proof the sufficient conditions.\\
Sept1,
for a fixed map $\beta\in GL(V)$. We define
$\rho:\frkg\longrightarrow\frkg\frkl(V)$ as follow
\begin{equation}\label{t1eq0}
d^sv(x)=\beta^{1+s}\rho(x)\beta^{-2-s}v,\forall v\in V,x\in\frkg.
\end{equation}
By straightforward computations, we have
$$(\bar{\alpha}\circ d^sv(x)=d^sv(\alpha(x))=\beta^{1+s}\rho(\alpha(x))\beta^{-2-s}v,$$
$$d^{s+1}\circ\bar{\alpha}(v)(x)=d^{s+1}v(x)=\beta^{2+s}\rho(x)\beta^{-3-s}v.$$
according to iii), we have:
\begin{equation}\label{t1eq1}
\rho(\alpha(x))\circ\beta=\beta\circ\rho(x).
\end{equation}
Sept2, for $\forall x,y\in\frkg, \eta\in C^1(\frkg;V)$, we define
$[\cdot,\cdot]:\frkg\wedge\frkg\longrightarrow\frkg$ by
\begin{equation}\label{t1eq2}
\langle \eta,[x,y]\rangle =\beta^{2+s}\rho(x)\beta^{-3-s}\eta(\alpha(y))-\beta^{2+s}\rho(y)\beta^{-3-s}\eta(\alpha(x))-d^s\eta(x,y).
\end{equation}
When $\xi\in\frkg^*$, according to (\ref{t1eq2}), we have
\begin{equation}\label{t1eq3}
 \langle \xi,[x,y]\rangle=-d\xi(x,y).
\end{equation}
According to (\ref{t1eq2}),(\ref{t1eq1})and iii), we have
\begin{eqnarray*}
\langle \eta,[\alpha(x),\alpha(y)]\rangle&=&\beta^{2+s}\rho(\alpha(x))\beta^{-3-s}\eta(\alpha^2(y))
-\beta^{2+s}\rho(\alpha(y))\beta^{-3-s}\eta(\alpha^2(x))-d^s\eta(\alpha(x),\alpha(y))\\
&&=\beta^{3+s}\rho(x)\beta^{-4-s}\eta(\alpha^2(y))
-\beta^{3+s}\rho(y)\beta^{-4-s}\eta(\alpha^2(x))-\bar{\alpha}\circ d^s\eta(x,y)\\
&&=\beta^{3+s}\rho(x)\beta^{-4-s}\bar{\alpha}(\eta)(\alpha(y))
-\beta^{3+s}\rho(y)\beta^{-4-s}\bar{\alpha}(\eta)(\alpha(x))-d^{s+1}\bar{\alpha}(\eta)(x,y)\\
&&=\langle \bar{\alpha}(\eta),[x,y]\rangle\\
&&=\langle \eta,\alpha([x,y])\rangle.
\end{eqnarray*}
So, we have
\begin{equation}\label{t1eq4}
 \alpha([x,y])=[\alpha(x),\alpha(y)].
\end{equation}
Sept3, for any $v\in C^0(\frkg;V)=V$, by i), (\ref{t1eq2}), (\ref{t1eq0}) and (\ref{t1eq1}), we have
\begin{eqnarray*}
0&=&d^s\circ d^sv(x,y)\\
&=&\beta^{2+s}\rho(x)\beta^{-3-s}d^sv(\alpha(y))-\beta^{2+s}\rho(y)\beta^{-3-s}d^sv(\alpha(x))-d^sv([x,y])\\
&=&\beta^{1+s}\rho(\alpha(x))\rho(y)\beta^{-3-s}v-\beta^{1+s}\rho(\alpha(y))\rho(x)\beta^{-3-s}v-\beta^{1+s}\rho([x,y])\beta^{-2-s}v.
\end{eqnarray*}
We get
\begin{equation}\label{t1eq5}
 \rho(\alpha(x))\rho(y)-\rho(\alpha(y))\rho(x)=\rho([x,y])\beta.
\end{equation}
Sept4, for any $\xi\in\frkg^*,\eta\in C^1(\frkg;V)$, according to ii), (\ref{t1eq3}) and (\ref{t1eq2}), we have
\begin{eqnarray*}
d^s(\xi\diamond\eta)(x,y,z)&=&d\xi\diamond\bar{\alpha}(\eta)(x,y,z)-\bar{\alpha}(\xi)\diamond d^{s+1}\eta(x,y,z)\\
&&=d\xi(x,y)\bar{\alpha}(\eta)(z)-d\xi(x,z)\bar{\alpha}(\eta)(y)+d\xi(y,z)\bar{\alpha}(\eta)(x)\\
&&-\bar{\alpha}(\xi)(x)d^{s+1}\eta(y,z)+\bar{\alpha}(\xi)(y)d^{s+1}\eta(x,z)-\bar{\alpha}(\xi)(z)d^{s+1}\eta(x,y)\\
&&=\beta^{3+s}\rho(x)\beta^{-4-s}(\xi\diamond\eta)(\alpha(y),\alpha(z))-\beta^{3+s}\rho(y)\beta^{-4-s}(\xi\diamond\eta)(\alpha(x),\alpha(z))\\
&&+\beta^{3+s}\rho(z)\beta^{-4-s}(\xi\diamond\eta)(\alpha(x),\alpha(y))-\xi\diamond\eta([x,y],\alpha(z))\\
&&+\xi\diamond\eta([x,z],\alpha(y))-\xi\diamond\eta([y,z],\alpha(x)).
\end{eqnarray*}
So, for any $\omega\in C^2(\frkg;V)$, we have
\begin{eqnarray*}
d^s\omega(x,y,z)&=&\beta^{3+s}\rho(x)\beta^{-4-s}\omega(\alpha(y),\alpha(z))-\beta^{3+s}\rho(y)\beta^{-4-s}\omega(\alpha(x),\alpha(z))\\
&&+\beta^{3+s}\rho(z)\beta^{-4-s}\omega(\alpha(x),\alpha(y))-\omega([x,y],\alpha(z))\\
&&+\omega([x,z],\alpha(y))-\omega([y,z],\alpha(x)).
\end{eqnarray*}
For any $\eta\in C^1(\frkg;V)$, according to i), we have
\begin{eqnarray*}
0&=&d^s\circ d^s\eta(x,y,z)\\
&=&\eta([[x,y],\alpha(z)]+[[y,z],\alpha(x)]+[[z,x],\alpha(y)])
\end{eqnarray*}
Then, we have
\begin{equation}\label{t1eq6}
 [[x,y],\alpha(z)]+[[y,z],\alpha(x)]+[[z,x],\alpha(y)]=0.
\end{equation}
So, according to (\ref{t1eq4}) and (\ref{t1eq6}), we have:
$(\frkg,[\cdot,\cdot],\alpha)$ is a  Hom-Lie algebra;\\
according to (\ref{t1eq1}) and (\ref{t1eq5}), we have: $\rho$ is a
representation of $(\frkg,[\cdot,\cdot],\alpha)$ on the vector space
$V$ with respect to $\beta$.\qed

\end{document}